\documentclass[11pt]{article}
\usepackage{amsmath}
\usepackage{amssymb, amscd, amsthm,amsfonts}
\usepackage[all]{xy}
\usepackage[dvips]{graphicx}
\usepackage{verbatim}
\usepackage[perpage,symbol*]{footmisc}
\usepackage[dvips]{color}
\newtheorem{theorem}{Theorem}
\newtheorem{lemma}{Lemma}

\newtheorem{corollary}{Corollary}
\newtheorem{conjecture}{Conjecture}

\begin{document}
	
\pagestyle{myheadings}

\title
{\bf The backtrackless circuits of $D(k,q)$ obtained by lifting\footnote{This work was supported by the National Natural Science Foundation of China (No. 61977056).}}

\author{Xiaoyan Cheng, Chen Zhang, Yuansheng Tang\footnote{Corresponding author.} \\
{\it\small School of Mathematical Sciences, Yangzhou University, Jiangsu, China}\\
and Ming Xu
{\it\small \footnote{Email addresses: xycheng@yzu.edu.cn(X. Cheng), cz12024@163.com(C. Zhang), ystang@yzu.edu.cn(Y. Tang), mxu@szcu.edu.cn(M.Xu).}}\\
{\it\small Suzhou City University,Jiangsu, China}}
\date{}

\maketitle	

{\noindent\small{\bf Abstract:}
For integer $k\geq2$ and power $q$ of prime $p$, the algebraic bipartite graph $D(k,q)$, proposed by Lazebnik and Ustimenko in 1993, is $q$-regular, edge-transitive and of girth at least $k+4$. For $q\geq4$, it was conjectured in 1995 that the exact girth of $D(k,q)$, denoted by $g(D(k,q))$, is equal to $g_k$, where $g_k$ is defined as $k+5$ if $2\nmid k$ and as $k+4$ if $2|k$. It is now known that this conjecture is true when $g_k$ is of form $hp^t$, where $h$ is a factor of $4(q-1)$. Sun proved in her thesis (2017) that, for $n\not\equiv 1(\bmod\ 4)$, any backtrackless circuit (BC) of length $l$ in $D(n,q)$ can be lifted to BCs of length $pl$ in $D(k,q)$, where $k$ is $p(n+5)-5$ if $n\equiv 3(\bmod\ 4)$ and $p(n+4)-4$ otherwise. This result implies that, provided $g(D(n,q))=g_n$, we must have $g(D(k,q))=g_k$ for any $k$ in the set $\{p^sg_n-5|s\geq 0\}$ if $n\equiv 3(\bmod\ 4)$ and for any $k$ in the set $\{p^sg_n-5|s\geq 0\}\cup\{p^sg_n-4|s\geq 0\}$ otherwise, where $D(1,q)$ is defined as the graph $D(2,q)$ if any. Sun's proof for this result is simplified in this paper, with the help of some new automorphisms of the graph $D(k,q)$ that are combinations of the known ones. A new proof for this result is also given in the last section of this paper by using a class of homogeneous polynomials.}
\vspace{1ex}


\section{Introduction}

In this paper we consider only simple graphs, which are undirected without loops and multiple edges.
For any graph $G$, its vertex set and edge set are denoted by $V(G)$ and $E(G)$, respectively.
For any vertex $v\in V(G)$, its {\it neighbourhood} $N(v)$ is the set of vertices which are adjacent to $v$.
   $G$ is said to be {\it edge-transitive} if, for any two edges $e_1$, $e_2$, that there exists an automorphism $\phi$ of $G$ such that $\phi$ maps the ends of $e_1$ into those of $e_2$.
   A {\it backtrackless walk} (BW) means a sequence $v_1v_2\cdots$ of vertices of $G$ such that any two consecutive vertices in the sequence are adjacent and $v_j\neq v_{j-2}$ for $j\geq 3$.
   For $n\geq 3$, a sequence $v_1v_2\cdots v_n$ of vertices is called a {\it backtrackless circuit} (BC) of length $n$
   if both $v_1v_2\cdots v_n$ and $v_3v_4\cdots v_nv_1v_2$ are BWs.
   The {\it girth} of a non-forest graph $G$, denoted by $g(G)$, is the length of the shortest
   BCs. Clearly, any BC in $G$ of length $g(G)$
   is a cycle since it must consist of distinct vertices.
   We call $G$ a {\it lift} (or {\it cover}) of another graph $G'$ if there exists a surjective mapping
$\theta: V(G)\rightarrow V (G')$, which satisfies the two conditions:
(i) $\theta$ preserves adjacencies, i.e., $\theta(u)\theta(v)\in E(G')$ whenever $uv\in E(G)$;
(ii) For any vertex $v\in V(G)$, the restriction of $\theta$ to $N(v)$ is a bijection between $N(v)$ and $N(\theta(v))$.
Clearly, the girth of a graph is not greater than the girth of its lift.

Graphs with large girth and a high degree of symmetry have been applied to variant problems in many fields, such as extremal graph theory, coding theory and cryptography. For $k\geq 2$ and prime power $q$, the algebraic bipartite graph $D(k,q)$, proposed by Lazebnik and Ustimenko in \cite{Lazebnik1993} and
investigated quite well in literature (c.f. \cite{Lazebnik93}-\cite{XCT22}), is $q$-regular, edge-transitive
and of girth at least
\begin{align*}
g_k=\left\{
\begin{array}{ll}
  k+4, & \text{if }k\text{ is even}, \\
  k+5, & \text{if }k\text{ is odd}.
\end{array}
\right.
\end{align*}
For the exact girth of $D(k,q)$, the following conjecture
   \begin{conjecture}\label{conmain}
     For $k\geq 2$ and $q\geq4$, the girth of $D(k,q)$ is equal to $g_k$.
   \end{conjecture}
   \noindent
was proposed in \cite{Fredi95}, and shown to be valid in \cite{Fredi95} for the case that $g_k/2$ is a factor of $q-1$, in \cite{XWY14} for the case that $g_k/2$ is a power of $p$, and in \cite{XWY16} for the case that $g_k/2$ is a product of a factor of $q-1$ and a power of $p$, respectively.
By modifying an idea from \cite{Fredi95}, it is proved in \cite{Sun17} that Conjecture~\ref{conmain} is valid
when $q$ is a power of odd prime $p$ and $g_k$ is a product of a factor of $4(q-1)$ and a power of $p$.
Recently, the same result was proved independently in \cite{XCT23} with a different method.
On the other hand, if the characteristic of $\mathbb{F}_q$ is $p$ and the graph $D(n,q)$ has girth $g_{n}$,
it is also proved in \cite{Sun17}, for any $s\geq 0$ and
\begin{align*}
k=\left\{
\begin{array}{ll}
  g_np^s-5, & \text{if }n\equiv 3(\bmod 4), \\
  g_np^s-4, & \text{otherwise},
\end{array}
\right.
\end{align*}
that the graph $D(k,q)$ has girth $g_k=g_np^s$ and some of its girth cycles can be constructed by lifting the girth cycles of $D(n,q)$.
Since it may be inconvenient for some readers to access \cite{Sun17}, in
the present paper we show for this result a proof modified slightly from that given in \cite{Sun17} and another proof
by using a compact expression of some BWs given in \cite{XCT22} and \cite{XCT23}.

Instead of working on the original form of the graph $D(k,q)$, we will deal with the bipartite graph
$\Lambda_{k,q}$, which is isomorphic to $D(k,q)$ and defined as follows (c.f. \cite{XWY14}, \cite{XWY16} and \cite{XCT22}).
The left part of vertices of $\Lambda_{k,q}$, denoted by $L_k$, consists of the vectors $[l]=(l_0,l_1,l_2,\ldots,l_k)\in\mathbb{F}_q^{k+1}$ with $l_1=l_2$.
The right part of vertices of $\Lambda_{k,q}$, denoted by $R_k$, consists of the vectors $\langle r \rangle=(r_0,r_1,r_2,\ldots,r_k)\in\mathbb{F}_q^k$ with $r_1=0$. Two vertices $(l_0,l_1,\ldots,l_k)\in L_k$ and $(r_0,r_1,\ldots,r_k)\in R_k$ are adjacent in $\Lambda_{k,q}$ if and only if, for $2\leq i \leq k$,
   \begin{align}\label{adj_con}
      l_i+r_i=\begin{cases}
         r_0l_{i-2} \text{ ~~if } i\equiv2,3 \mod 4,\\
         l_0r_{i-2} \text{ ~~if } i\equiv0,1 \mod 4.
      \end{cases}
   \end{align}
Clearly, for $2\leq n<k$ the graph $\Lambda_{k,q}$ is a lift of the graph $\Lambda_{n,q}$.

This paper is arranged as follows.
In Section~2, for each edge $e\in E(\Lambda_{k,q})$ we show an automorphism of $\Lambda_{k,q}$ that
   transforms the edge $[0]\langle 0\rangle$ into $e$, this implies that $\Lambda_{k,q}$ is edge-transitive,
   a well-known result. With the help of these automorphisms, we show the main result of this paper in Section~3,
   where some BCs in $\Lambda_{k,q}$ are constructed by lifting. In particular, for fixed prime power $q$ the correctness of Conjecture~\ref{conmain} for any given $k$ implies the same kind of correctness for infinite $k$'s.
   In Section~4, we give another proof for the main result of this paper by use of a compact expression for the vertices on some BWs of $\Lambda_{k,q}$.

\section{Automorphisms}

In this section, we consider to give some new automorphisms of $\Lambda_{k,q}$ by combining the known ones
(c.f. \cite{Lazebnik93}, \cite{Fredi95} and \cite{Lazebnik17}).
For any given edge $e$, one of the new automorphisms transforms the edge $[0]\langle 0\rangle$ into $e$.

At first, for $j\geq 0$, we define $\mu(j)$ as the following set of indices
\begin{align}\label{auto01}
  \mu(j)=\left\{
  \begin{array}{ll}
    \{i|0\leq i\leq j, i\equiv 0,1(\bmod 4)\}, & \text{if } j\equiv 0,1(\bmod\ 4),\\
    \{i|2\leq i\leq j, i\equiv 2,3(\bmod 4)\}\cup\{0\}, & \text{if } j\equiv 2,3(\bmod\ 4).
  \end{array}
  \right.
\end{align}
Clearly, for $l\geq 0$ we have
\begin{gather}
\label{at0}
\mu(l+2)=\left\{
\begin{array}{ll}
\{i+2:i\in\mu(l)\}\cup\{0\},&\text{if }l\equiv 0,1\bmod 4,\\
\{i+2:i\in\mu(l)\setminus\{0\}\}\cup\{0,1\},&\text{if }l\equiv 2,3\bmod 4,
\end{array}
\right.\\
\label{at1}
\mu(l)\cap(\mu(l)+2)=\left\{
\begin{array}{ll}
\emptyset,&\text{if }l\equiv 0,1\bmod 4,\\
\{2\},&\text{if }l\equiv 2,3\bmod 4,
\end{array}
\right.\\
\label{at2}
\mu(l)\cup(\mu(l)+2)=\left\{
\begin{array}{ll}
\{0,1,\ldots,l+2\}\setminus\{l+1\},&\text{if }l\equiv 0\bmod 4,\\
\{0,1,\ldots,l+2\},&\text{if }l\equiv 1\bmod 4,\\
\{0,1,\ldots,l+2\}\setminus\{1,l+1\},&\text{if }l\equiv 2\bmod 4,\\
\{0,1,\ldots,l+2\}\setminus\{1\},&\text{if }l\equiv 3\bmod 4.
\end{array}
\right.
\end{gather}
For any edge $e=[l^*]\langle r^*\rangle\in L_k\times R_k$ of $\Lambda_{k,q}$, let $\theta_{e}$ denote the map defined on $L_k\cup R_k$ with $\theta_{e}([l])=[\hat{l}]\in L_k$ and $\theta_{e}(\langle r\rangle)=\langle \hat{r}\rangle\in R_k$ such that
\begin{align}\label{auto02}
\hat{l}_j=l^*_j+\sum_{0\leq i\leq j}f_{j,i}l_i,\
\hat{r}_j=r^*_j+\sum_{0\leq i\leq j}g_{j,i}r_i,
\end{align}
where
\begin{align}\label{auto00'}
f_{1,0}=f_{2,0}=r^*_0, \ \ g_{1,0}=f_{2,1}=0
\end{align}
and, for $j\geq 0$,
\begin{align}
&f_{j,j}=g_{j,j}=1,\label{auto00}\\
&f_{j,i}=0 \text{ for }i\not\in\mu(j) \text{ if }j\equiv 0,1(\bmod\ 4), \\
&g_{j,i}=0 \text{ for }i\not\in\mu(j) \text{ if }j\equiv 2,3(\bmod\ 4).\label{m0}
\end{align}

From (\ref{auto02}), (\ref{auto00}), $\hat{r}_1=r_1^*+g_{1,0}r_0+g_{1,1}r_1=0$ and
\begin{gather*}
\hat{l}_{2}=f_{2,0}l_0+f_{2,1}l_1+f_{2,2}l_2
=r^*_0l_0+l_2=r^*_0l_0+l_1=f_{1,0}l_0+f_{1,1}l_1=\hat{l}_1,
\end{gather*}
we see that the map $\theta_{e}$ is indeed a bijection on $L_k$ and $R_k$, respectively.
To enable the map $\theta_{e}$ preserving the adjacency, we manage to impose some additional conditions on it.

Let $[l]\langle r\rangle \in L_k\times R_k$ be an arbitrary edge of $\Lambda_{k,q}$.
For $j\equiv 2,3(\bmod\ 4)$, according to (\ref{adj_con}), (\ref{auto01}), (\ref{auto02}) and (\ref{auto00})-(\ref{m0}) we see
\begin{gather*}
\Big(r^*_j+\sum_{i\in\mu(j)}g_{j,i}r_i\Big)+\Big(l^*_j+\sum_{0\leq i\leq j}f_{j,i}l_i\Big)
=\hat{r}_j+\hat{l}_j
\end{gather*}
should equal to
\begin{align*}
&\hat{r}_0\hat{l}_{j-2}=(r^*_0+r_0)\Big(l^*_{j-2}+\sum_{i\in\mu(j-2)}f_{j-2,i}l_i\Big)\\
=&(r^*_0+r_0)l^*_{j-2}+\sum_{i\in\mu(j-2)}f_{j-2,i}(r^*_0+r_0) l_i\\
=&l^*_{j}+r^*_j+l^*_{j-2}r_0+\sum_{i\in\mu(j-2)}f_{j-2,i}(r^*_0l_i+l_{i+2}+r_{i+2})\\
=&\Big(r^*_j+l^*_{j-2}r_0+\sum_{i\in\mu(j-2)}f_{j-2,i}r_{i+2}\Big)+\Big(l^*_{j}+
\sum_{i\in\mu(j-2)}f_{j-2,i}(r^*_0l_i+l_{i+2})\Big),
\end{align*}
and thus according to (\ref{at0})-(\ref{at2}) we should set, for $t\geq 0$ and $\tau\in\{0,1\}$, that
\begin{align}\label{auto06}
g_{4t+\tau+2,i}=\left\{
\begin{array}{ll}
  l^*_{4t+\tau}, & \text{if }i=0, \\
  f_{4t+\tau,i-2}, & \text{otherwise},
\end{array}
\right.
\end{align}
\begin{align}\label{auto07}
f_{4t+\tau+2,i}=\left\{
\begin{array}{ll}
  r^*_{0}f_{4t+\tau,i}, & \text{if }i\in\mu(4t+\tau), \\
  f_{4t+\tau,i-2}, & \text{if }i\in 2+\mu(4t+\tau),\\
  0, & \text{otherwise}.
\end{array}
\right.
\end{align}
We note that, for $j=4t+\tau+2=2$, (\ref{auto07}) coincides those already defined by (\ref{auto00'}) and (\ref{auto00}).
For $j\geq 4$ and $j\equiv 0,1(\bmod\ 4)$, according to (\ref{adj_con}), (\ref{auto01}), (\ref{auto02}) and (\ref{auto00})-(\ref{m0}) we see
\begin{gather*}
\Big(l^*_j+\sum_{i\in\mu(j)}f_{j,i}l_i\Big)+\Big(r^*_j+\sum_{0\leq i\leq j}g_{j,i}r_i\Big)=\hat{l}_j+\hat{r}_j
\end{gather*}
should equal to
\begin{align*}
&\hat{l}_0\hat{r}_{j-2}=(l^*_0+l_0)\Big(r^*_{j-2}+\sum_{i\in\mu(j-2)}g_{j-2,i}r_i\Big)\\
=&(l^*_0+l_0)r^*_{j-2}+\sum_{i\in\mu(j-2)}g_{j-2,i}(l^*_0+l_0) r_i\\
=&l^*_{j}+r^*_j+r^*_{j-2}l_0+\sum_{i\in\mu(j-2)}g_{j-2,i}(l^*_0r_i+l_{i+2}+r_{i+2})\\
=&\Big(l^*_j+r^*_{j-2}l_0+\sum_{i\in\mu(j-2)}g_{j-2,i}l_{i+2}\Big)+\Big(r^*_{j}+
\sum_{i\in\mu(j-2)}g_{j-2,i}(l^*_0r_i+r_{i+2})\Big)
\end{align*}
and thus according to (\ref{at0})-(\ref{at2}) we should set, for $t\geq 1$ and $\tau\in\{0,1\}$, that
\begin{align}\label{auto08}
f_{4t+\tau,i}=\left\{
\begin{array}{ll}
  r^*_{4t+\tau-2}, & \text{if }i=0, \\
  g_{4t+\tau-2,0}, & \text{if }i=1, \\
  g_{4t+\tau-2,i-2}, & \text{otherwise},
\end{array}
\right.
\end{align}
\begin{align}\label{auto09}
g_{4t+\tau,i}=\left\{
\begin{array}{ll}
  l^*_0g_{4t+\tau-2,2}+g_{4t+\tau-2,0}, & \text{if }i=2,\\
  l^*_{0}g_{4t+\tau-2,i}, & \text{if }i\in\mu(4t+\tau-2)\setminus\{2\}, \\
  g_{4t+\tau-2,i-2}, & \text{if }i\in \mu(4t+\tau)\setminus\{0,1\},\\
  0, & \text{otherwise}.
\end{array}
\right.
\end{align}

Hereafter, we assume that $e=[l^*]\langle r^*\rangle \in L_k\times R_k$ is an edge of $\Lambda_{k,q}$ and the map $\theta_{e}$ satisfies (\ref{auto02})-(\ref{auto09}).
Then one can deduce easily the following theorem, which implies that $\Lambda_{k,q}$ is edge-transitive.
\begin{theorem}\label{lem300}
The map $\theta_{e}$
is an automorphism of $\Lambda_{k,q}$ that transforms the edge $[0]\langle 0\rangle$ into $e$.
\end{theorem}

Since the automorphism $\theta_e$ defined above is of iterative manner, in the end of this section  we give it directly by the entries
of the two ends of $e$.
From (\ref{auto06}) and (\ref{auto08}), we see easily the following lemma.
\begin{lemma}\label{lem301}
Let $t\geq 0$ and $\tau\in\{0,1\}$. Then, $f_{4t+1,4t}=g_{4t+3,4t+2}=r^*_0$, $g_{4t+2+\tau,0}=l^*_{4t+\tau}$
 and, for $0\leq s<t$,
\begin{gather*}
  f_{4t+\tau,4s}=g_{4t+2+\tau,4s+2}=r^*_{4t-2+\tau-4s},\\
  f_{4t+\tau,4s+1}=g_{4t+2+\tau,4s+3}=l^*_{4t-4+\tau-4s}.
\end{gather*}
\end{lemma}
\noindent
Furthermore, according to Lemma~\ref{lem301}, from (\ref{auto07}) and (\ref{auto09}) one can deduce easily the following two lemmas, respectively.
\begin{lemma}\label{lem601}
Let $t\geq 0$. Then, $f_{4t+2,4t+1}=0$, $f_{4t+2,4t}=f_{4t+3,4t+1}=f_{4t+3,4t+2}=r^*_0$, $f_{4t+3,4t}=(r^*_0)^2$
 and, for $0\leq s<t$ and $\tau\in\{0,1\}$,
\begin{gather*}
  f_{4t+2+\tau,4s}=r^*_0r^*_{4t-2+\tau-4s},\\
  f_{4t+2+\tau,4s+1}=r^*_0l^*_{4t-4+\tau-4s},\\
  f_{4t+2+\tau,4s+2}=r^*_{4t-2+\tau-4s},\\
  f_{4t+2+\tau,4s+3}=l^*_{4t-4+\tau-4s}.
\end{gather*}
\end{lemma}
\begin{lemma}\label{lem602}
Let $t\geq 0$, $\tau\in\{0,1\}$, $r^*_{-2}=1$ and $r^*_{-1}=r^*_0$. Then, $g_{4t+4,4t+3}=g_{4t+4+\tau,1}=0$,
$g_{4t+5,4t+3}=l^*_0$, $g_{4t+5,4t+4}=r^*_0$,
$g_{4t+4+\tau,0}=l^*_0l^*_{4t+\tau}$, $g_{4t+4+\tau,2}=l^*_0r^*_{4t-2+\tau}+l^*_{4t+\tau}$
and, for $0\leq s<t$,
\begin{gather*}
  g_{4t+4+\tau,4s+3}=l^*_0l^*_{4t-4+\tau-4s},\\
  g_{4t+4+\tau,4s+4}=r^*_{4t-2+\tau-4s},\\
  g_{4t+4+\tau,4s+5}=l^*_{4t-4+\tau-4s},\\
  9_{4t+4+\tau,4s+6}=l^*_0r^*_{4t-6+\tau-4s}.
\end{gather*}
\end{lemma}

Clearly, the automorphism $\theta_e$ of $\Lambda_{k,q}$ can also be defined by (\ref{auto02})-(\ref{m0}) and Lemmas~\ref{lem301}, \ref{lem601}, \ref{lem602}.

\section{BCs Constructed By Lifting}

In this section, we will investigate the BCs of $\Lambda_{k,q}$ obtained by lifting, with the help of the automorphisms introduced in the last section.

For $2\leq n\leq k$, let $L_{k,n}$ ($R_{k,n}$, respectively) denote the set of tuples in $L_k$ ($R_k$, respectively) whose
first $n+1$ entries are 0.

Now we assume that $e=[l^*]\langle r^*\rangle\in L_{k,n}\times R_{k,n}$ is a given edge of $\Lambda_{k,q}$
and $\phi_e$ is the automorphism defined in the last section. Clearly, for any $[l]\in L_{k,n}$ and $\langle r\rangle\in R_{k,n}$, we have $\phi_e([l])\in L_{k,n}$ and $\phi_e(\langle r\rangle)\in R_{k,n}$.
We note that Lemmas~\ref{lem301}, \ref{lem601} and \ref{lem602} are used frequently in the proof of the following lemma without notification.

\begin{lemma}\label{lemv00}
For $j\geq 0$, let
\begin{align}\label{auto20}
\delta_{j,n}=\left\{
\begin{array}{ll}
  n+3, & \text{if } n\equiv \tau+1(\bmod\ 4),\\
  n+5, & \text{if }n\equiv \tau+2(\bmod\ 4),\\
  n+4, & \text{otherwise},\\
\end{array}
\right.
\end{align}
where $\tau\in\{0,1\}$ is the parity number of $j$, i.e., $\tau$ is 0 if $j$ is even and 1 otherwise.
Then,
\begin{align}\label{auto30}
f_{j,i}=g_{j,i}=0\text{ for }0<j-i<\nabla
\end{align}
is valid for $\nabla=\delta_{j,n}$.
\end{lemma}
\begin{proof}
According to $l^*_t=r^*_t=0$ for $0\leq t\leq n$,
one can check easily that (\ref{auto30}) is valid for $\nabla=n+3$.
Now we manage to enlarge such kind of upper bound $\nabla$.

If $n\equiv\tau+1(\bmod\ 4)$, then we have $f_{j,j-n-3}=g_{j,j-n-3}=r^*_{n+1}$ and thus the bound $\nabla$ cannot be replaced with a larger number for this case.

If $n\equiv\tau(\bmod\ 4)$, then we have $f_{j,j-n-3}=g_{j,j-n-3}=l^*_n$, $f_{j,j-n-4}=g_{j,j-n-4}=r^*_{n+2}$.
Hence, we see that (\ref{auto30}) is valid for $\nabla=n+4$ if $n\equiv\tau(\bmod\ 4)$.

Assume $n\equiv\tau+3(\bmod\ 4)$.
If $j\equiv0,1(\bmod\ 4)$, then $f_{j,j-n-3}=0$, $g_{j,j-n-3}=l^*_0r^*_{n-1}$, $f_{j,j-n-4}=g_{j,j-n-4}=l^*_{n+1}$. If $j\equiv2,3(\bmod\ 4)$, then $f_{j,j-n-3}=r^*_0r^*_{n-1}$, $g_{j,j-n-3}=0$, $f_{j,j-n-4}=g_{j,j-n-4}=l^*_{n+1}$.
Hence, we see that (\ref{auto30}) is valid for $\nabla=n+4$ if $n\equiv\tau+3(\bmod\ 4)$.

Assume $n\equiv\tau+2(\bmod\ 4)$.
If $j\equiv0,1(\bmod\ 4)$, then $f_{j,j-n-3}=0$, $g_{j,j-n-3}=l^*_0l^*_{n-2}$, $f_{j,j-n-4}=0$, $g_{j,j-n-4}=l^*_0r^*_{n}$, $f_{j,j-n-5}=g_{j,j-n-5}=l^*_{n+2}$. If $j\equiv2,3(\bmod\ 4)$, then $f_{j,j-n-3}=r^*_0l^*_{n-2}$, $g_{j,j-n-3}=0$, $f_{j,j-n-4}=r^*_0r^*_{n}$, $g_{j,j-n-4}=0$, $f_{j,j-n-5}=g_{j,j-n-5}=l^*_{n+2}$.
Hence, we see that (\ref{auto30}) is valid for $\nabla=n+5$ if $n\equiv\tau+2(\bmod\ 4)$.

The proof is complete.
\end{proof}

For positive integers $n,s$ with $n\not\equiv 1(\bmod\ 4)$, let
\begin{align}
\xi_{n,s}=\left\{
\begin{array}{ll}
  (n+4)(s+1)-4, & \text{if }2|n, \\
  (n+5)(s+1)-5, & \text{if }n\equiv 3(\bmod\ 4).
\end{array}
\right.
\end{align}
The case $n\equiv 1(\bmod\ 4)$ is excluded because of $\Lambda_{4i+2,q}\cong q\Lambda_{4i+1,q}$ (cf. \cite{Lazebnik96}, \cite{Lazebnik04}). Then, we have the following lemma.

\begin{lemma}\label{lemv03}
For any $s\geq 1$ and $\xi_{n,s}\geq j_0>j_1>\cdots>j_s>n$, we have
\begin{align}
f_{j_0,j_1}\cdots f_{j_{s-1},j_s}=0.
\end{align}
\end{lemma}
\begin{proof}
Let $\xi_{n,0}=n$. Clearly, there must be some positive integers $k,t$ with $k,t\leq s$ such that $\xi_{n,t}\geq j_{k-1}>j_{k}>\xi_{n,t-1}$.
We note that
\begin{align*}
\xi_{n,t}=\left\{
\begin{array}{cc}
(t+1)n+4t,&\text{if }n\equiv0(\bmod\ 4),\\
(t+1)n+4t,&\text{if }n\equiv2(\bmod\ 4),\\
(t+1)n+5t,&\text{if }n\equiv3(\bmod\ 4),\\
\end{array}
\right.
\end{align*}
\begin{align*}
\xi_{n,t}-\xi_{n,t-1}=\left\{
\begin{array}{cc}
n+4,&\text{if }n\equiv0(\bmod\ 4),\\
n+4,&\text{if }n\equiv2(\bmod\ 4),\\
n+5,&\text{if }n\equiv3(\bmod\ 4).\\
\end{array}
\right.
\end{align*}
If $j_{k-1}\equiv\xi_{n,t}(\bmod\ 2)$, then from
\begin{align*}
\delta_{j_{k-1},n}=\left\{
\begin{array}{cc}
n+4,&\text{if }n\equiv0(\bmod\ 4),\\
n+5,&\text{if }n\equiv2(\bmod\ 4),\\
n+5,&\text{if }n\equiv3(\bmod\ 4),\\
\end{array}
\right.
\end{align*}
we have $0<j_{k-1}-j_k\leq \xi_{n,t}-(\xi_{n,t-1}+1)<\delta_{j_{k-1},n}$ and thus from Lemma~\ref{lemv00} we see $f_{j_{k-1},j_k}=0$.
If $j_{k-1}\equiv\xi_{n,t}-1(\bmod\ 2)$, then from
\begin{align*}
\delta_{j_{k-1},n}=\left\{
\begin{array}{cc}
n+4,&\text{if }n\equiv0(\bmod\ 4),\\
n+3,&\text{if }n\equiv2(\bmod\ 4),\\
n+4,&\text{if }n\equiv3(\bmod\ 4),\\
\end{array}
\right.
\end{align*}
we have $0<j_{k-1}-j_k\leq (\xi_{n,t}-1)-(\xi_{n,t-1}+1)<\delta_{j_{k-1},n}$ and thus from Lemma~\ref{lemv00} we also see $f_{j_{k-1},j_k}=0$.
\end{proof}

Now we can prove the main result of this paper.
\begin{theorem}\label{thm004}
Assume $n$ is a positive integer with $n\not\equiv 1(\bmod\ 4)$.
If $\Lambda_{n,q}$ has a BC of length $2g$, then for $k=\xi_{n,p-1}$ the graph $\Lambda_{k,q}$ has a BC of length $2pg$, where
$p$ is the characteristic of $\mathbb{F}_q$.
\end{theorem}
\begin{proof}
Let $\Gamma=[l^{(1)}]\langle r^{(1)}\rangle[l^{(2)}]\langle r^{(2)}\rangle\cdots[l^{(g)}]\langle r^{(g)}\rangle$ be a lifting in $\Lambda_{k,q}$ of a BC of length $2g$ in $\Lambda_{n,q}$.
Since the graph $\Lambda_{k,q}$ is edge-transitive, without loss of generality we assume further that
$[l^{(1)}]=[0]$ and $\langle r^{(1)}\rangle=\langle 0\rangle$ are the all zero-vectors.
Let $e=[l^*]\langle r^*\rangle\in L_k\times R_k$ be the edge in $\Lambda_{k,q}$ with $l^*_0=r^*_0=0$ and
$[l^*]\langle r^{(g)}\rangle\in E(\Lambda_{k,q})$.
Then, we must have $e\in L_{k,n}\times R_{k,n}$.

Since, for each $j_0>n$ and $s\geq 1$, from (\ref{auto02}) and (\ref{auto00}) we can show easily by induction that the $(j_0+1)$-th entry of $\theta_e^{s+1}([0])=\theta_e^{s}([l^*])$ is
\begin{align*}
(s+1)l^*_{j_0}+\sum_{l=1}^{s}\binom{s+1}{l+1}\sum_{j_0>j_1>\cdots>j_l\geq 0}f_{j_0,j_1}\cdots f_{j_{l-1},j_l}l^*_{j_l},
\end{align*}
then, from $[l^*]\in L_{k,n}$, Lemma~\ref{lemv03} and $p|\binom{p}{l+1}$ for $1\leq l<p-1$,
we have $\theta_e^{p}([0])=[0]$, and therefore the concatenation of
$$\Gamma,\theta_e(\Gamma),\cdots,\theta_e^{p-1}(\Gamma)$$
forms a BC in $\Lambda_{k,q}$ of length $2pg$.
\end{proof}

From Theorem~\ref{thm004} and the bound $g(\Lambda_{k,q})\geq g_k$, one can deduce the following corollary simply.
\begin{corollary}\label{cor001}
If $g(\Lambda_{4t-5,q})=4t$, then $g(\Lambda_{4tp^s-5,q})=4tp^s$ for any $s\geq 0$.
If $g(\Lambda_{2t-4,q})=2t$, then $g(\Lambda_{2tp^s-5,q})=g(\Lambda_{2tp^s-4,q})=2tp^s$ for any $s\geq 0$.
\end{corollary}

We note that the previous proof for Theorem \ref{thm004} is indeed a simple modification of that given in \cite{Sun17}.
In the next section, we will give another proof for it.

\section{Another Proof of the Main Result}
In this section, we consider to give another proof for Theorem~\ref{thm004} by using
the compact expression for the vertices on some BWs of $\Lambda_{k,q}$ given in \cite{XCT22} and \cite{XCT23}.

At first, we introduce the homogeneous polynomials defined in \cite{XWY14}.
For indeterminates $w_1,\ldots,w_n$, let
\begin{align*}
\rho_0(w_1,\ldots,w_n)=w_1\cdots w_n
\end{align*}
	and, for $1 \leq s\leq \lfloor\frac{n}{2}\rfloor$, let $\rho_s(w_1,\ldots,w_n)$ denote the homogeneous polynomial of order $n-2s$ defined by
\begin{align}
\rho_s(w_1,\ldots,w_n)=\sum_{1\leq i_1<\cdots<i_s\leq n-s}\frac{\prod_{j=1}^n w_j}{\prod_{j=1}^s w_{i_j+j-1} w_{i_j+j}},
\end{align}
	where each term in the summation is a product of the remaining elements in the sequence $w_1,\ldots,w_n$ after deleting from it $s$ disjoint pairs $\{w_i,w_{i+1}\}$ of consecutive elements.
	If $n<2s$ or $s<0$, $\rho_s(w_1,\ldots,w_n)$ is defined as 0. For the null sequence $\eta$, we define $\rho_s(\eta)$ as
	\begin{align*}
		\rho_s(\eta)=\begin{cases}
			1 \text{ ~~if } s=0,\\
			0 \text{ ~~if } s\neq 0.
		\end{cases}
	\end{align*}
Furthermore, for any integers $i,t$ with $i\geq 1$ and $2i$-tuple $\alpha=(w_1,\ldots,w_{2i})$, we define notations
\begin{align}
  \Delta^{\alpha}_{2t-1,1}=& \rho_{i-t}(w_2,w_3,\ldots,w_{2i}), \label{delta01}\\
  \Delta^{\alpha}_{2t-1,0}=& \rho_{i-t}(w_1,w_2,\ldots,w_{2i-1}), \label{delta02}\\
  \Delta^{\alpha}_{2t,1}=& \rho_{i-1-t}(w_2,w_3,\ldots,w_{2i-1}), \label{delta03}\\
  \Delta^{\alpha}_{2t,0}=& \rho_{i-t}(w_1,w_2,\ldots,w_{2i}). \label{delta04}
\end{align}
We note that, for each nonzero homogenous polynomial in the right-hand sides of (\ref{delta01})-(\ref{delta04}), its order is equal to the first subscript in the corresponding notation in the left-hand sides.

Let $\Gamma=[l^{(1)}]\langle r^{(1)}\rangle[l^{(2)}]\langle r^{(2)}\rangle\cdots$ be a given BW of $\Lambda_{k,q}$ with $[l^{(1)}]=[0]$ and $\langle r^{(1)}\rangle=\langle 0\rangle$. For $i\geq1$, let $x_i$ and $y_i$ denote the first entries (or colors) of $[l^{(i)}]$ and $\langle r^{(i)}\rangle$, respectively, and write
	\begin{align}\label{relationxy}
		u_i=x_{i+1}-x_{i}, ~~~~~~~v_i=y_{i+1}-y_{i}.
	\end{align}
Clearly, we have $u_i\neq0$ and $v_i\neq0$.
The following lemma (c.f. \cite{XCT22} and \cite{XCT23}) shows a compact expression for the vertices on $\Gamma$.
	\begin{lemma}\label{path}
		Assume $i\geq 1$.
Let $l_j^{(i+1)}$ be the $(j+1)$-th entry of $[l^{(i+1)}]$ for $j\geq  0$. Then, for $s\geq 0$ we have
		\begin{align}
			l_{4s}^{(i+1)}&=\Delta^{\epsilon_i}_{2s+1,0},\label{4j}\\
			l_{4s+1}^{(i+1)}&=\Delta^{\epsilon_i}_{2s+2,1},\label{4j+1}\\
			l_{4s+2}^{(i+1)}&=y_{i+1}l_{4s}^{(i+1)}-\Delta^{\epsilon_i}_{2s+2,0},\label{4j+2}\\
			l_{4s+3}^{(i+1)}&=y_{i+1}l_{4s+1}^{(i+1)}-\Delta^{\epsilon_i}_{2s+3,1},\label{4j+3}
		\end{align}
where $\epsilon_i$ is the $2i$-tuple $(u_1,v_1,\ldots,u_i,v_i)$.
	\end{lemma}

If the first $2i$ vertices on $\Gamma$ form a BC in $\Lambda_{k,q}$ of length $2i$, we also call it a BC of {\it type} $\epsilon_i=(u_1,v_1,\ldots,u_i,v_i)$.
We note that
     \begin{align}\label{4j-1}
        \Delta^{\epsilon_i}_{1,1}=v_1+\cdots+v_i=y_{i+1}=y_1=0
     \end{align}
is always a necessary condition for $\Lambda_{k,q}$ to have a BC of type $\epsilon_i$.
Since each BC in $\Lambda_{4t+1,q}$ can be lifted naturally to some BCs in $\Lambda_{4t+2,q}$ of the same type (cf. \cite{Lazebnik96}, \cite{Lazebnik04}, \cite{XCT23}), from Lemma \ref{path} one can deduce simply the following theorem.

\begin{theorem}
\label{thm002}
Let $\alpha=(\omega_1,\ldots,\omega_{2n})$ be a $2n$-tuple over $\mathbb{F}_q^*$.
\begin{enumerate}
\item For $1\leq s\leq \lfloor(n-2)/2\rfloor$, the graph $\Lambda_{4s,q}$ has a BC of type $\alpha$ if and only if $\Delta^{\alpha}_{i,\sigma}=0$ for any $\sigma\in\{0,1\}$ and $1\leq i\leq 2s+1$.
\item For $1\leq s\leq \lfloor(n-2)/2\rfloor$, the graph $\Lambda_{4s-1,q}$ has a BC of type $\alpha$ if and only if $\Delta^{\alpha}_{2s+1,1}=0$ and $\Delta^{\alpha}_{i,\sigma}=0$ for any $\sigma\in\{0,1\}$ and $1\leq i\leq 2s$.
\item For $k\in\{4s-2,4s-3\}$ with $1\leq s\leq \lfloor(n-1)/2\rfloor$, the graph $\Lambda_{k,q}$ has a BC of type $\alpha$ if and only if $\Delta^{\alpha}_{i,\sigma}=0$ for any $\sigma\in\{0,1\}$ and $1\leq i\leq 2s$.
\end{enumerate}
\end{theorem}

Now we investigate further the properties of the homogenous polynomials defined above. At first, we note for $0\leq j\leq n$ that one can show easily (c.f. \cite{XCT22} and \cite{XCT23}) that
   \begin{align}
\rho_{n-j}(\omega_1,\ldots,\omega_{2n})&=\sum_{1\leq s_1\leq t_1<s_2\leq t_2<\cdots < s_j\leq t_j\leq n}
         \prod_{k=1}^j\omega_{2s_k-1}\omega_{2t_k},\label{0b0}\\
\rho_{n-j}(\omega_1,\ldots,\omega_{2n+1})&=\sum_{1\leq s_0\leq t_1< s_1\leq t_2<s_2\leq \cdots \leq t_j< s_j\leq n+1}
         \omega_{2s_0-1}\prod_{k=1}^j\omega_{2t_k}\omega_{2s_k-1}.\label{0b1}
   \end{align}
Then, one can deduce the following lemma easily.

\begin{lemma}\label{lem400}
For any $\sigma_0\in\{0,1\}$, positive integer $t$ and tuples $\alpha$, $\beta$ of even length,
we have
\begin{align}\label{cc00}
  \Delta^{\alpha\beta}_{t,\,\sigma_0}&
  =\sum_{t_1+t_2=t}\Delta^{\alpha}_{t_1,\,\sigma_0}\Delta^{\beta}_{t_2,\,\sigma_0+t_1\bmod 2},
\end{align}
where $\alpha\beta$ denotes the concatenation of $\alpha$ and $\beta$.
\end{lemma}
\begin{proof}
We prove (\ref{cc00}) only for the case $\sigma_0=0$ and $2|t$, the other cases can be proved similarly.

Assume $\alpha=(a_1,a_2,\ldots,a_{2m})$ and $\beta=(b_1,b_2,\ldots,b_{2n})$.
Then, according to (\ref{0b0}) and (\ref{0b1}) we have
\begin{align*}
&\Delta^{\alpha\beta}_{2l,0}=\rho_{m+n-l}(a_1,\ldots,a_{2m},b_1,\ldots,b_{2n})\\
=&\sum_{i+j=l}\sum_{1\leq s_1\leq t_1<\cdots < s_i\leq t_i\leq m}\prod_{k=1}^ia_{2s_k-1}a_{2t_k}\sum_{1\leq s'_1\leq t'_1<\cdots < s'_j\leq t'_j\leq n}\prod_{k=1}^jb_{2s'_k-1}b_{2t'_k}\\
&\ \ \ \ \ \ \ \ +\sum_{i+j=l-1}\sum_{1\leq s_0\leq t_1< s_1\leq \cdots \leq t_i< s_i\leq m}a_{2s_0-1}
\prod_{k=1}^ia_{2t_k}a_{2s_k-1}\\
&\ \ \ \ \ \ \ \ \ \ \ \ \ \ \ \ \ \ \ \ \ \ \ \ \ \ \ \ \ \
\sum_{1\leq t'_0< s'_1\leq t'_1< \cdots < s'_j\leq t'_j\leq n}b_{2t'_0}\prod_{k=1}^jb_{2s'_k-1}b_{2t'_k}\\
=&\sum_{i+j=l}\rho_{m-i}(a_1,\ldots,a_{2m})\rho_{n-j}(b_1,\ldots,b_{2n})\\
&\ \ \ \ \ \ \ \ \ \ \ \ \ \ \ \ +\sum_{i+j=l-1}\rho_{m-1-i}(a_1,\ldots,a_{2m-1})\rho_{n-1-j}(b_2,\ldots,b_{2n})\\
=&\sum_{i+j=l}\Delta^{\alpha}_{2i,0}\Delta^{\beta}_{2j,0}+\sum_{i+j=l-1}\Delta^{\alpha}_{2i+1,0}\Delta^{\beta}_{2j+1,1}\\
=&\sum_{t_1+t_2=2l}\Delta^{\alpha}_{t_1,\,0}\Delta^{\beta}_{t_2,\,t_1\bmod 2},
\end{align*}
i.e. (\ref{cc00}) is valid when $\sigma_0=0$ and $2|t$.
\end{proof}

For positive integers $j,t$, let $\Omega_{j}(t)$ denote the set of $j$-tuple $\tau=(\tau_1,\cdots,\tau_j)$ of nonnegative integers with $\tau_1+\cdots+\tau_j=t$. We note that $|\Omega_j(t)|=\binom{j+t-1}{t}$.
Let $\Omega^*_{j}(t)$ denote the set of tuples in $\Omega_{j}(t)$ consisting of positive integers.
For any tuple $\tau\in\Omega_j(t)$, let $\tau^*$ denote the tuple obtained from $\tau$ by deleting all of the zero entries.
For $\zeta\in\Omega_s^*(t)$, let
$\Theta_j(\zeta)$ be the set of tuples $\tau\in\Omega_j(t)$ with $\tau^*=\zeta$.
Clearly, we have the following partition
\begin{align}\label{cc502}
\Omega_{j}(t)=\bigcup_{1\leq s\leq j, \zeta\in \Omega^*_{s}(t)}\Theta_j(\zeta).
\end{align}

Since, for any $\zeta\in\Omega_s^*(t)$ with $1\leq s\leq j$, $\Theta_j(\zeta)$ is the set of tuples obtained by inserting $j-s$ zero entries into $\zeta$, we have
\begin{align}
|\Theta_j(\zeta)|=|\Omega_{s+1}(j-s)|=\binom{(s+1)+(j-s)-1}{j-s}=\binom{j}{s}.
\end{align}
For $\tau=(\tau_1,\cdots,\tau_j)\in\Omega_j(t)$ and $\sigma_0\in\{0,1\}$, let $(\sigma_1,\ldots,\sigma_j)$ denote the $j$-tuple over $\{0,1\}$ defined by
\begin{align}\label{cc503}
\sigma_i=\sigma_i(\tau,\sigma_0)\equiv\sigma_{i-1}+\tau_{i-1}\equiv\sigma_0+\sum_{1\leq k<i}\tau_i(\bmod\ 2),\ i=1,2,\ldots,j,\end{align}
where $\tau_0=0$.
\begin{lemma}
\label{lem002}
Let $\tau'$ be a tuple obtained from
  $\tau\in\Omega_j(t)$ by deleting an even entry.
  Then, for any $\sigma_0\in\{0,1\}$, we have
\begin{align*}
\sigma_i(\tau',\sigma_0)=\left\{
\begin{array}{ll}
\sigma_{i}(\tau,\sigma_0),&\text{ if }1\leq i<k,\\
\sigma_{i+1}(\tau,\sigma_0),&\text{ if }k\leq i<j,
\end{array}
\right.
\end{align*}
where $k$ is the index of the deleted entry of $\tau$.
\end{lemma}
\begin{proof}
Clearly, we have
$\sigma_i(\tau,\sigma_0)=\sigma_i(\tau',\sigma_0)$ for $1\leq i<k$.
For $l>k$, from $2|\tau_k$ we see that
\begin{align*}
\sigma_{l}(\tau,\sigma_0)\equiv\sigma_0+\sum_{1\leq i<l}\tau_i
\equiv\sigma_0+\sum_{1\leq i<l,i\neq k}\tau_i\equiv\sigma_0+\sum_{1\leq i<l-1}\tau'_i(\bmod\ 2)
\end{align*}
is equal to $\sigma_{l-1}(\tau',\sigma_0)$.
\end{proof}

Assume that $\alpha$ is a tuple of even length. For $j\geq 1$, let $\alpha^j$ denote the concatenation of $j$ copies of $\alpha$.
For any $\sigma_0\in\{0,1\}$ and positive integer $t$,
according to Lemmas~\ref{lem400}, \ref{lem002} and (\ref{cc502})-(\ref{cc503}), one can deduce easily that
\begin{align}
\Delta^{\alpha^j}_{t,\sigma_0}&=
\sum_{\tau\in\Omega_j(t)}\prod_{i=1}^j\Delta^{\alpha}_{{\tau}_i,\sigma_i(\tau,\sigma_0)}\nonumber\\
&=\sum_{1\leq s\leq j}\sum_{\zeta\in\Omega^*_s(t)}
\sum_{\tau\in\Theta_j(\zeta)}\prod_{i=1}^j\Delta^{\alpha}_{{\tau}_i,\sigma_i(\tau,\sigma_0)}\nonumber\\
&=\sum_{1\leq s\leq j}\binom{j}{s}\sum_{\zeta\in\Omega^*_s(t)}
\prod_{i=1}^s\Delta^{\alpha}_{{\zeta}_i,\sigma_i(\zeta,\sigma_0)},\label{cc501}
\end{align}
where the fact $\Delta^{\alpha}_{0,\sigma}=1$ is used in the last equality.

Now we can prove the following theorem, an alternate of Theorem~\ref{thm004}.
\begin{theorem}
Assume that $n$ is a positive integer with $n\not\equiv 1(\bmod\ 4)$ and
$\Lambda_{n,q}$ has a BC of type $\alpha$. Then, for $k=\xi_{n,p-1}$, the graph $\Lambda_{k,q}$ has a BC of type $\alpha^p$, where
$p$ is the characteristic of $\mathbb{F}_q$.
\end{theorem}
\begin{proof}
At first, we assume that $n=2l$ is even.
According to Theorem~\ref{thm002} we have $\Delta^{\alpha}_{i,\sigma_0}=0$ for $1\leq i\leq l+1$ and $\sigma_0\in\{0,1\}$.
Hence, from (\ref{cc501}) we see $\Delta^{\alpha^p}_{i,\sigma_0}=0$ for $1\leq i\leq p(l+2)-1$
and $\sigma_0\in\{0,1\}$. Therefore, from Theorem~\ref{thm002} and $$(k+2)/2=(\xi_{2l,p-1}+2)/2=(p(2l+4)-4+2)/2=p(l+2)-1$$ we see that $\Lambda_{k,q}$ has a BC of type $\alpha^p$.

Now we assume $n=4s-1$.
According to Theorem~\ref{thm002} we have $\Delta^{\alpha}_{2s+1,1}=0$ and $\Delta^{\alpha}_{i,\sigma_0}=0$ for $1\leq i\leq 2s$ and $\sigma_0\in\{0,1\}$.

Let $\zeta\in\Omega^*_p(t)$ be an arbitrary tuple with $\prod_{i=1}^p\Delta^{\alpha}_{{\zeta}_i,\sigma_i(\zeta,\sigma_0)}\neq 0$.
Clearly, each of the entries of $\zeta$ is not smaller than $2s+1$.
Let
$i_1<i_2<\cdots<i_a$ be the indices of the entries which are equal to $2s+1$.
Then, we have $\sigma_{i_b}(\zeta,\sigma_0)=0$ for $1\leq b\leq a$.
From (\ref{cc503}) one can deduce easily that there must be some indices $j_1,j_2,\ldots,j_{a-1}$ such that $i_b<j_b<i_{b+1}$ and $\zeta_{j_b}\geq 2s+3$ for $1\leq b<a$.
Therefore, we have
$$t=\zeta_1+\zeta_2+\cdots+\zeta_p\geq 2(s+1)p-1.$$
Furthermore, if $\sigma_0=1$, from (\ref{cc503}) we see that there must be some index $j_0$ such that $1\leq j_0<i_1$ and $\zeta_{j_0}\geq 2s+3$. Hence, we have $t\geq 2(s+1)p$ if $\sigma_0=1$.

Therefore, from (\ref{cc501}) we see $\Delta^{\alpha^p}_{2(s+1)p-1,1}=0$ and
$\Delta^{\alpha^p}_{i,\sigma_0}=0$ for $1\leq i\leq 2(s+1)p-2$
and $\sigma_0\in\{0,1\}$.
Hence, from Theorem~\ref{thm002} and
$$(k+1)/2=(\xi_{4s-1,p-1}+1)/2=(p(4s-1+5)-5+1)/2=2(s+1)p-2$$
we see that
$\Lambda_{k,q}$ has a BC of type $\alpha^p$.
\end{proof}

Clearly, Corollary \ref{cor001} can also be deduced simply from this theorem.

\bigskip
\noindent
{\bf Acknowledgements}: The authors wish to thank Professor Felix Lazebnik for providing the preprint of
\cite{Sun17}.

\end{document}